\documentclass[11pt]{amsart}

\usepackage{amsmath,amsthm, amscd, amssymb, amsfonts}
\usepackage{epsfig}
\input xy
\xyoption{all}

\newcommand{\Mo}{{\mathcal M}}

\newcommand{\Ss}{{\mathcal S}}
\newcommand{\ot}{{\otimes}}
\newcommand{\otb}{{\overline{\otimes}}}

\newcommand{\Ac}{{\mathcal A}}

\newcommand{\ele}{{\mathcal L}}

\newcommand{\ca}{{\mathcal C}}

\newcommand{\ku}{{\Bbbk}}

\newcommand{\Na}{{\mathbb N}}

\newcommand{\uno}{{\mathbf 1}}

\newcommand{\id}{\mbox{\rm id\,}}

\newcommand{\Res}{\mbox{\rm Res\,}}
\newcommand{\Ind}{\mbox{\rm Ind\,}}

\newcommand\Rep{\operatorname{Rep}}

\newcommand\co{\operatorname{co}}

\newcommand\Hom{\operatorname{Hom}}
\newcommand\uhom{\underline{\Hom}}

\newcommand{\End}{\operatorname{End}}
\newcommand\uend{\underline{\End}}

\newcommand\St{\operatorname{Stab}}

\newcommand{\Id}{\mathop{\rm Id}}

\renewcommand{\_}[1]{\mbox{$_{\left( #1 \right)}$}}

\theoremstyle{plain}

\numberwithin{equation}{section}

\newtheorem{teo}{Theorem}[section]

\newtheorem{lema}[teo]{Lemma}

\newtheorem{cor}[teo]{Corollary}

\newtheorem{prop}[teo]{Proposition}

\theoremstyle{definition}

\newtheorem{defi}[teo]{Definition}

\theoremstyle{remark}

\newtheorem{rmk}[teo]{Remark}

\def\pf{\begin{proof}}

\def\epf{\end{proof}}

\theoremstyle{remark}

\begin{document}

\title[Dynamical twists in  Hopf algebras]
{Dynamical twists in  Hopf algebras}
\author[ Mombelli]{
Juan Mart\'\i n Mombelli }
\thanks{This work was partially supported by
 CONICET, Secyt (UNC)}
\address{Facultad de Matem\'atica, Astronom\'\i a y F\'\i sica, Universidad Nacional de
\newline
\indent C\'ordoba, CIEM -- CONICET.
\newline \indent Medina Allende s/n, (5000) Ciudad Universitaria, C\'ordoba, Argentina}
\email{mombelli@mate.uncor.edu \newline \indent\emph{URL:}\/
http://www.mate.uncor.edu/mombelli}
\begin{abstract} We establish a bijective correspondence between
gauge equi\-valence classes of dynamical twists in a
finite-dimensional Hopf algebra $H$ based on a finite abelian
group $A$ and equivalence classes of pairs $(K,
\{V_{\lambda}\}_{\lambda\in \widehat{A}})$, where $K$ is an
$H$-simple left $H$-comodule semisimple algebra and
$\{V_{\lambda}\}_{\lambda\in \widehat{A}}$ is a family of
irreducible representations satisfying certain conditions. Our
results generalize the results obtained by Etingof-Nikshych on the
classification of dynamical twists in group algebras.

\end{abstract}

\date{\today}
\maketitle

\setcounter{tocdepth}{1}
\section{Introduction}

The notion of dynamical twist introduced in \cite{B}, see also
\cite{BBB}, is a generalization of Drinfeld's notion of twist to
the dynamical setting. More precisely, if $A$ is a finite Abelian
subgroup of the group of group-like elements of a Hopf algebra $H$
a dynamical twist for the pair $(H,A)$ is a function
$J:\widehat{A} \to H\ot H$ satisfying certain equations. If $A$ is
trivial then a dynamical twist is just a usual twist. In
\cite{EN2} for any dynamical twist the authors endowed the product
$H\ot_{\ku} \End_{\ku}(A)$ with a nontrivial weak Hopf algebra
structure. One of the main properties that they prove is that if
$H$ is quasitriangular with $R$-matrix $R$ then
$\mathcal{R}(\lambda)=J^{-1}(\lambda)^{21} RJ(\lambda)$ satisfies
the dynamical quantum Yang-Baxter equation. See \cite{E} for a
comprehensive presentation of the dynamical quantum Yang-Baxter
equation.

\medbreak

Etingof and Nikshych \cite{EN1} classify dynamical twists for
group algebras of finite groups in terms of group data. In this
paper we are concerned with the classification of dynamical twists
over any finite-dimensional Hopf algebras. The approach of this
work owns a lot to \cite{EN1}, however there are some differences;
the langauge of module categories is used with profit, and we make
use of the stabilizers for Hopf algebra actions introduced by M.
Yan and Y. Zhu \cite{YZ}.

\medbreak

The paper is organized as follows. In subsection \ref{notation} we
introduce the basic notation and conventions, also the main tools
that will be needed further. We recall the definition of modules
over a tensor category and some results from \cite{AM} on modules
over the category of representation of a Hopf algebra. We briefly
explain the stabilizers for Hopf algebra actions and some of their
properties.

\smallbreak

In section \ref{ddd} we begin with the definition of dynamical
twists over a Hopf algebra and gauge equivalence between them. In
subsection \ref{modcat-twist} from a dynamical twist over a Hopf
algebra $H$ we construct a module category over $\Rep(H)$. This
module category will be important to understand dynamical twists.

\smallbreak

In subsection \ref{datum} following the spirit of the work
\cite{EN1} we introduce the definition of dynamical datum. The
main new ingredient that appear is the Yan-Zhu stabilizer. In
subsection \ref{exchange} we show how to construct a dynamical
twists from a dynamical datum in such a way that equivalence
classes of dynamical data give the same gauge equivalence of
dynamical twists. A converse of this result is proved in
Proposition \ref{dynamdat}. Finally in subsection \ref{main0} we
prove that the constructions explained above are one the inverse
of each other. This is our main result stated in Theorem
\ref{main}.

\smallbreak

In section \ref{examples} we compute some examples. Specifically,
in subsection \ref{k=a} we compute the dynamical twist in the case
when $K$ is the group algebra of the group $A$, and in subsection
\ref{taft-ex} we show a one-parameter family of dynamical twists
for the Taft Hopf algebras.


\subsection*{Acknowledgments} The author is grateful to Nicol\'as
Andruskiewitsch for many remarks that improve the presentation of
the paper.

\section{Preliminaries}
\subsection{Notation}\label{notation} Throughout this paper $\ku$ will denote an
algebraically closed field of characteristic $0$, and $H$ shall
denote a finite-dimensional Hopf algebra over $\ku$ with counit
$\varepsilon$ and antipode $\Ss$. We shall denote by $\Rep(H)$ the
tensor category of finite-dimensional left $H$-modules.

\medbreak

If $A$ is a finite Abelian group the group of characters of $A$
will be denoted by $\widehat{A}$. For any $\lambda\in\widehat{A}$
we sometimes denote by $\ku_{\lambda}$ the one-dimensional
representation afforded by $\lambda$.

\medbreak

Let $K$ be a left $H$-comodule algebra with coaction $\delta:K\to
H\ot K$. An $H$-\emph{costable ideal} of $K$ is a two-sided ideal
$I$ of $K$ such that $\delta(I)\subseteq H\ot I$. We shall say
that $K$ is \emph{$H$-simple} if there exists no non-trivial
$H$-costable ideal of $K$.

We denote by ${}^{H\!}\Mo_{K}$ the category of left $H$-comodules,
right $K$-modules $M$ such that the left $K$-module structure $
M\ot_{\ku} K\to M$ is an $H$-comodule map. If $S$ is another left
$H$-comodule algebra, the category $^{H\!}_{S\!}\Mo_K$ will denote
the category of $(S,K)$-bimodules $M$ with a left $H$-comodule
structure such that both actions are morphisms of $H$-comodules.

\begin{lema}\label{comodd}  Let $P \in {}^{H\!}\Mo_K$ .
Then $\End_K(P)$ is a left $H$-comodule algebra via
$\delta:\End_K(P)\to H\ot_{\ku} \End_K(P) $, $T\mapsto T\_{-1}\ot
T\_0$, determined by
\begin{equation}\label{h-comod} \langle\alpha, T\_{-1}\rangle\,
T_0(p)=\langle\alpha, T(p\_0)\_{-1}\Ss^{-1}(p\_{-1})\rangle\,
T(p\_0)\_0,\end{equation}  $T\in\End_K(P),$ $p\in P$, $\alpha\in
H^*$. Moreover $P \in ^{H\!}_{S\!}\Mo_K$, where $S=\End_K(P)$.
\end{lema}
\pf See \cite[Lemma 1.26]{AM}.\epf

We shall need later the following Frobenius reciprocity.

\begin{lema}\label{frobenius} Let $R$ be a subalgebra of $H$. For
every left $R$-module $W$ and a left $H$-module $V$ there is a
natural isomorphism
$$\Hom_R(W, \Res^H_R V)\simeq \Hom_H(\Ind^H_R  W, V).  $$
\end{lema}
\pf See, for example, \cite[Lemma 3.1]{AN}. \epf

\subsection{Module categories}\label{subsection-modcat} We briefly recall the definition of module
category and the definition introduced by Etingof-Ostrik of
\emph{exact} module categories. We refer to \cite{O1}, \cite{O2},
\cite{eo}.

Let us fix $\ca$ a finite tensor category. A \emph{module
category} over $\ca$ is a collection $(\Mo,\otb, m, l)$ where
$\Mo$ is an Abelian category, $\otb:\ca\times\Mo\to \Mo$ is an
exact bifunctor, associativity and unit isomorphisms
$m_{X,Y,M}:(X\ot Y)\otb M\to X\otb ( Y\otb M)$, $l_M:\uno\otb M\to
M$, $X,Y\in \ca$, $M\in\Mo$, such that
\begin{align}\label{modcat1} m_{X,Y,Z\otb M}\, m_{X\ot
Y,Z,M}&=(\id_X\ot\, m_{Y,Z,M})\,m_{X,Y\ot Z,M}\,(a_{X,Y,Z}\ot\id_M),\\
\label{modcat2} (\id_X\ot\, l_M)\, m_{X,\uno,M}&= r_X\ot\id_M,
\end{align}
for all $X,Y, Z\in \ca$, $M\in\Mo$. Sometimes we shall simply say
that $\Mo$ is a module category omitting to mention $\otb, m $ and
$l$ whenever no confusion arises.

\medbreak

In this paper we further assume that all module categories have
finitely many isomorphism classes of simple objects.

\medbreak

Let $\Mo$, $\Mo'$ be two module categories over $\ca$. A module
functor between $\Mo$ and $\Mo'$ is a pair $(F,c)$ where $F:\Mo\to
\Mo'$ is a functor and $c_{X,M}:F(X\otb M)\to X\otb F(M)$ are
natural isomorphisms such that
\begin{align}\label{modfun1} m'_{X,Y,F(M)}\, c_{X\ot Y,M}&=
(\id_X\ot\,
c_{Y,M})\, c_{X,Y\otb M} F(m_{X,Y,M}),\\
\label{modfun2} l'_{F(M)}\, c_{\uno, M} &= F(l_M),
\end{align}
for all $X,Y\in \ca$, $M\in\Mo$. Two module categories $\Mo$ and
$\Mo'$ are \emph{equivalent} if there exists a module functor
$(F,c)$ where $F$ is an equivalence of categories. The module
structure over $\Mo\oplus \Mo'$ is defined in an obvious way. A
module category is \emph{indecomposable} if it is not equivalent
to the direct sum of two non-trivial module categories.

A module category $\Mo$ is \emph{exact} \cite{eo} if for any
projective object $P\in \ca$ and any $M\in \Mo$ the object $P\otb
M$ is again projective.

\medbreak

We recall the definition of the internal Hom. This object is an
important tool in the study of module categories. See for example
\cite{O1}, see also \cite{eo}.

\smallbreak

Let $\Mo$ be an exact module category over $\ca$. Let $M_1, M_2
\in \Mo$. The functor $X\mapsto \Hom_{\Mo}(X\otb M_1, M_2)$ is
representable and an object $\uhom(M_1, M_2)$ representing this
functor is called the \emph{internal Hom} of $M_1$ and $M_2$. See
\cite{eo, O1} for details. Thus
$$
\Hom_{\Mo}(X\otb M_1, M_2) \simeq \Hom_{\ca}(X, \uhom(M_1, M_2))
$$
for any $X\in \ca$, $M_1, M_2 \in \Mo$.

\smallbreak   The internal Hom $\uhom(M,M)$ of an object $M\in
\Mo$ is an algebra in $\ca$. The multiplication is constructed as
follows. Denote by $$ev_M:\uend(M)\ot M\to M$$ the
\emph{evaluation map} \label{evaluation} obtained as the image of
the identity under the isomorphism
$$\Hom_{\ca}(\uend(M),\uend(M))\simeq  \Hom_{\Mo}(\uend(M)\ot M,M).$$
Thus the product $\mu:\uend(M)\ot \uend(M)\to \uend(M)$ is defined
as the image of the map
\begin{equation}\label{mult-intt} ev_{M}
(\id\ot\, ev_{M})\, m_{\uend(M),\uend(M),M}
\end{equation}
under the isomorphism
\begin{align*}\Hom_{\ca}(\uend(M)\ot\uend(M),
\uend(M)) \simeq\Hom_{\Mo}((\uend(M)\ot\uend(M)) \ot M, M).
\end{align*}
For details we refer to \cite{O1}.

\subsection{Yan-Zhu stabilizers} In \cite{YZ} the authors introduce
a notion of stabilizers for Hopf algebra actions generalizing the
existing notion for groups. In \cite{AM} these objects and a
slight extension thereof, called \emph{Yan-Zhu stabilizers}, were
used in connection with module categories over Hopf algebras. We
recall the definition and some of its properties.

\medbreak

Let $K$ be a finite-dimensional left $H$-comodule algebra an $V,$
$W$ two left $K$-modules. The \emph{Yan-Zhu stabilizer}
$\St_K(V,W)$ is defined as the intersection
$$ \St_K(V,W)= \Hom_K(H^* \ot V,H^* \ot W)
\cap \ele\big(H^* \ot \Hom (V, W)\big).$$ Here the map $\ele:H^*
\ot \Hom (V, W)$ $\to$ $ \Hom(H^* \ot V,H^* \ot W)$ is defined by
$\ele(\gamma\ot T)(\xi\ot v)=\gamma\xi\ot T(v)$, for every
$\gamma, \xi\in H^*$, $T\in \Hom (V, W)$, $v\in V$.

\medbreak

The $K$-action on $H^*\ot V$ is given by
$$ k\cdot (\gamma\ot v)= k_{(-1)} \rightharpoondown \gamma\, \ot\,
k_{(0)}\cdot v,$$ for all $k\in K$, $\gamma\in H^*$, $v\in V$.
Here $\rightharpoondown:H\ot H^*\to H^*$ is the action defined by
$\langle h\rightharpoondown \gamma, t\rangle=\langle \gamma,
\Ss^{-1}(h)t\rangle$, for all $h, t\in H$, $\gamma\in H^*$. Also,
we denote $\St_K(V)=\St_K(V,V)$. The Yan-Zhu stabilizers can also
be presented as the coinvariants of a certain Hopf module, see
\cite[Prop. 2.10]{AM}.

\begin{prop}\cite[Prop. 2.7, Prop. 2.16]{AM}\label{stab-properties} The following assertions holds.

1. For any left $K$-module $V$ the stabilizer $\St_K(V)$ is a left
$H$-module algebra.

2. If $K$ is $H$-simple then
\begin{equation}\label{dim-stab} \dim(K) \dim(\St_K(V, W))= \dim(V)
\dim(W) \dim(H).
\end{equation}

3. For any $X\in \Rep(H)$ there are natural isomorphisms
$$\Hom_H(X,\St_K(V, W))\simeq \Hom_K(X\ot_{\ku} V, W).$$\qed
\end{prop}

In general the Yan-Zhu stabilizer is not easy object to compute,
however we have the following result. Let $\widetilde{H}\subseteq
H$ be a Hopf subalgebra and $K\supseteq R=K^{\co H}$ be a
Hopf-Galois extension over $\widetilde{H}$.

\begin{prop}\label{hopf-galois} For any $V, W$ left $K$-modules
there is an $H$-module isomorphism $\St_K(V,W)\simeq
\Hom_{\widetilde{H}}(H, \Hom_R(V, W)).$
\end{prop}
\pf See \cite[Theorem 2.23]{AM}.\epf

Let $K$ be a $H$-simple left $H$-comodule algebra with coaction
given by $\delta:K\to H\ot_{\ku} K$. The category of finite
dimensional left $K$-modules ${}_K\Mo$ is an exact module category
over $\Rep(H)$, see \cite[Prop. 1.20 (i)]{AM}. The action
$\otb:\Rep(H)\times {}_K\Mo\to {}_K\Mo$ is $X\otb V=X\ot_{\ku} V$,
where the left $K$-action on $X\ot_{\ku} V$ is given by $\delta$,
that is
\begin{equation}\label{mod-action} k\cdot (x\ot v)= k\_{-1}\cdot
x\ot\, k\_0\cdot v,
\end{equation}
for all $k\in K$, $x\in X$, $v\in V$. The associativity and unit
isomorphisms are canonical.

Proposition \ref{stab-properties} (3) implies that the Yan-Zhu
stabilizers are the internal Hom's of the module category
$_{K}\Mo$. As a consequence of this observation and  \cite[Th.
3.17]{eo} we have the following result.

\begin{cor}\label{modk} If $K$ is $H$-simple and $V$ is a left $K$-module,
then ${}_K\Mo\simeq \Rep(H)_{\St_K(V)}$ as module categories over
$\Rep(H)$.
\end{cor}
\pf See \cite[Corollary 3.2]{AM}. \epf

Let us recall a result from \cite{AM} that will be very useful
later.

\begin{teo}\cite[Theorem 3.3]{AM}\label{mod-overhopf} If $\Mo$ is an exact idecomposable
module category over $\Rep(H)$ then $\Mo\simeq {}_K\Mo$ for some
$H$-simple left comodule algebra $K$ with $K^{\co H}=\ku$.\qed
\end{teo}

Let $S$ be another $H$-simple left $H$-comodule algebra.
\begin{prop}\label{eq1}\cite[Prop. 1.24]{AM} The module categories ${}_K\Mo$, ${}_S\Mo$
over $\Rep(H)$ are equivalent if and only if there exists $P\in
{}^{H\!}\Mo_{S}$ such that $K\simeq \End_S(P)$ as $H$-module
algebras. Moreover the equivalence is given by $F:{}_K\Mo\to
{}_S\Mo$, $F(V)=P\ot_S V$, for all $V\in {}_K\Mo$. \qed
\end{prop}

\section{Dynamical twists for Hopf algebras}\label{ddd}

Hereafter $H$ will denote a finite-dimensional Hopf algebra and
$A\subseteq G(H)$ a finite Abelian subgroup of the group of
group-like elements of $H$. We first recall the definition of
dynamical twist over $H$, see \cite{B}, \cite{BBB}, \cite{EN1},
\cite{EN2}.

\begin{defi} Let $J:\widehat{A}\to H\ot H$ be a linear map with
invertible values. $J$ is a \emph{dynamical twist for} $H$ if for
any $\lambda\in \widehat{A}$ and $a\in A$

\begin{align}\label{d1} J(\lambda)(a\ot a)&=(a\ot a)J(\lambda),\\
\label{d2} \sum_{\mu\in \widehat{A}} (\Delta\ot\id)J(\lambda)\,
\left(J(\lambda\mu^{-1})\ot
P_{\mu}\right)&=(\id\ot \Delta)J(\lambda)\,(1\ot J(\lambda)),\\
\label{d3} (\varepsilon\ot\id)J(\lambda)&=(\id\ot
\varepsilon)J(\lambda)=1.
\end{align}
Here, for any $\mu \in \widehat{A}$,
$P_{\mu}=\displaystyle{\frac{1}{|A|} \sum_{a\in A} \mu(a^{-1})\,
a}\in \ku A$, is the minimal idempotent corresponding to the
character $\mu$.
\end{defi}


We will often use the notation: $J(\lambda)=J^1(\lambda)\ot\,
J^2(\lambda)=j^1(\lambda)\ot j^2(\lambda),$
$J(\lambda)^{-1}=J^{-1}(\lambda)\ot \,J^{-2}(\lambda),$
$\lambda\in\widehat{A}$.

\begin{defi} Two dynamical twists $J, J':\widehat{A}\to H\ot H$
are \emph{gauge equivalent} if there exists a map
$t:\widehat{A}\to H,$ called the gauge equivalence, with
invertible values such that for all $a\in A$, $\lambda \in
\widehat{A}$
\begin{align}\label{gauge3} \varepsilon(t(\lambda))&=1, \\
\label{gauge1} t(\lambda)\,a\,&=\, a\, t(\lambda), \\
\label{gauge2} J'(\lambda)&=\Delta(t(\lambda)^{-1})\, J(\lambda)
\, \sum_{\mu\in \widehat{A}} \left(t(\lambda\mu^{-1})\ot\, P_{\mu}
t(\lambda)\right).
\end{align}
\end{defi}
For any $\lambda\in \widehat{A}$ and $V$ a left $A$-module we
denote by $V[\lambda]$ the isotypic component of type $\lambda$,
that is
$$V[\lambda]=\{v\in V: a\cdot v=\lambda(a)\,v\,
\text{ for all } a\in A \}. $$ In particular if $X\in \Rep(H)$,
then $X[\lambda]$ makes sense by restriction of the action to
$\ku[A]$.

\medbreak

Now we will define an $H$-module algebra that will be useful
later. For any $\mu\in \widehat{A}$ define $$C(\mu)=\{f\in H^*:
f(ah)=\mu(a) f(h) \text{ for all } a\in A \}.$$ In another words
$C(\mu)= H^*[\lambda]$. We shall consider the left action of $H$
on $C(\mu)$ defined as follows. For $t, h\in H$ and $f\in C(\mu)$,
$(h\cdot f)(t)=f(th)$.

\medbreak

If $J:\widehat{A}\to H\ot H$ is a map that takes invertible
elements and satisfies \eqref{d1}, for each $\lambda\in
\widehat{A}$ we will define a product in $C(\varepsilon)$ as
follows. Let $f, g\in C(\varepsilon)$ then define
\begin{align}\label{product1} (f\cdot
g)(t)=f(J^{-1}(\lambda)h\_1)\; \,g(J^{-2}(\lambda)h\_2),
\end{align}
for all $h\in H$. Since $J$ commutes with $\Delta(a)$ for all
$a\in A$ then $f\cdot g\in C(\varepsilon)$. We shall denote by
$B_{\lambda}$ the space $C(\varepsilon)$ with this product.

\begin{lema}\label{homint11} If $J:\widehat{A}\to H\ot H$ is a
dynamical twist then $B_{\lambda}$ is an $H$-module algebra.
\end{lema}
\pf We only prove the associativity of the product. Let $f, g,
j\in B_{\lambda}$ and $h\in H$ then
\begin{align*} (f\cdot (g \cdot l))(h)&=f(J^{-1}(\lambda)h\_1)\;
(g \cdot l)(J^{-2}(\lambda)h\_2)\\
&=f(J^{-1}(\lambda)h\_1)\;
g(j^{-1}(\lambda)J^{-2}(\lambda)\_1h\_2)\;
l(j^{-2}(\lambda)J^{-2}(\lambda)\_2h\_3)\\
&= f(j^{-1}(\lambda)J^{-1}(\lambda)\_1h\_1)\,g(j^{-2}(\lambda)
J^{-1}(\lambda)\_2h\_2)\,l(J^{-2}(\lambda)h\_3)\\
&=((f\cdot g)\cdot l)(h).
\end{align*}
The third equality by \eqref{d2} since $l\in C(\varepsilon)$ and
$\varepsilon(P_{\mu})=\delta_{\mu,\varepsilon}$ for all $\mu\in
\widehat{A}$. \ \epf

A variation of the algebra $B_{\lambda}$ was considered in
\cite{EN1}. For usual twists the algebra $B_{\lambda}$ was first
considered in \cite{Mov} for the classification of twists in group
algebras. See also \cite{AEGN} and references therein.

\begin{rmk}\label{isod} There is an isomorphism of $H$-modules
$C(\lambda^{-1})\simeq(\Ind^H_A \lambda)^*$ given by
$\xi:C(\lambda^{-1})\to(\Ind^H_A \lambda)^*$,
$\xi(f)(\overline{h\ot 1})= f(\Ss(h))$, for all $h\in H$. Here
$\overline{h\ot 1}$ denotes the class of $h\ot 1$ in $H\ot_A
\ku_{\lambda}$.
\end{rmk}

\subsection{Module categories coming from dynamical twists}\label{modcat-twist}
For any dynamical twist we construct a semisimple module category
over $\Rep(H)$. The idea of relate dynamical twists with module
categories is due to Ostrik, see \cite[section 4.4]{O1}. In
\emph{loc. cit.} the author interpret the classification of
dynamical twists over group algebras given in \cite{EN1} as a
particular case of his classification of module categories over
group algebras.

\medbreak

Let $J:\widehat{A}\to H\ot H$ be a dynamical twist. Denote by
$\Mo^{(J)}$ the Abelian category of left $\ku[A]$-modules with the
following module category structure. Define
$\overline{\ot}:\Rep(H)\times\Mo^{(J)}\to \Mo^{(J)}$,
$X\overline{\ot} V:=X\ot_{\ku} V$, $X, Y\in \Rep(H)$, $V\in
\Mo^{(J)}$, where the $A$-module structure over $X\ot_{\ku} V$ is
given by the diagonal map. Define also $ m_{X,Y,V}:(X\ot Y)\ot
V\to X\otb (Y\otb V) $ by
\begin{align}\label{assoc-mc} m_{X,Y,V}(x\ot y\ot n)= J^{-1}(\lambda)\cdot
x\ot\, J^{-2}(\lambda)\cdot y\ot\,  n,\end{align} for any $x\in X,
y\in Y, n\in V[\lambda^{-1}]$. Since $J$ commutes with elements in
$A$ then $m_{X,Y,V}$ is an $\ku[A]$-module map.

\begin{lema}\label{mod-twist1} $(\Mo^{(J)}, \otb, m, \id)$ is an indecomposable
module category over $\Rep(H)$.
\end{lema}
\begin{proof} Let $X, Y, Z\in \Rep(H),$ $M\in \Mo^{(J)}$ and
$x\in X, y\in Y, z\in Z[\mu], n\in M[\lambda^{-1}]$. Then
\begin{align*} m_{X,Y,Z\ot M}\, &m_{X\ot Y,Z,M}(((x\ot y)\ot z)\ot
n)=\\&=m_{X,Y,Z\ot M}(J^{-1}(\lambda)\cdot (x\ot y)\ot
J^{-2}(\lambda)\cdot z\ot n)\\
&=j^{-1}(\lambda\mu^{-1})J^{-1}(\lambda)\_1\cdot x\ot
j^{-2}(\lambda\mu^{-1})J^{-1}(\lambda)\_2 \ot J^{-2}(\lambda)\cdot
z\ot n.
\end{align*}
Next,
\begin{align*}(\id_X\ot\, m_{Y,Z,M})&\,m_{X,Y\ot Z,M}((x\ot (y\ot z))\ot
n)=\\
&=(\id_X\ot\, m_{Y,Z,M})(J^{-1}(\lambda)\cdot x\ot
J^{-2}(\lambda) \cdot (y\ot z) \ot n\\
&=J^{-1}(\lambda) \cdot x\ot
j^{-1}(\lambda)J^{-2}(\lambda)\_1\cdot y\ot
j^{-2}(\lambda)J^{-2}(\lambda)\_2 \cdot z\ot n.
\end{align*}
Thus equation \eqref{d2} implies \eqref{modcat1}. Equation
\eqref{modcat2} follows immediately from \eqref{d3}. The
indecomposability of the module category  $\Mo^{(J)}$ follows as
in \cite[Prop. 1.18]{AM}.
\end{proof}

\begin{lema}\label{mod-twist2} If $J, \widetilde{J}$ are gauge equivalent
dynamical twists, then $\Mo^{(J)}\simeq \Mo^{(\widetilde{J})}$ as
module categories over $\Rep(H)$.
\end{lema}
\pf Let $t:\widehat{A}\to H$ be the gauge equivalence for $J$ and
$J'$. Define the module functor $(F,c):\Mo^{(J)}\to
\Mo^{(\widetilde{J})}$ as follows. For any $X\in \Rep(H), M\in
\Mo^{(J)}$, $F(M)=M$ and $c_{X,M}:X\ot_{\ku} M\to X\ot_{\ku} M$ is
defined by
$$c_{X,M}(x\ot m)= t(\lambda)\cdot x\ot
\, m,$$ for all $x\in X, m\in M[\lambda^{-1}]$. Since $t(\lambda)$
commutes with elements in $A$ the map $c_{X,M}$ is an $A$-module
morphism. Equation \eqref{modfun1} follows from \eqref{gauge2} and
equation \eqref{modfun2} follows from \eqref{gauge3}. Clearly
$(F,c)$ is an equivalence of module categories. \epf

The internal Hom of the module category $\Mo^{(J)}$ can be
explicitly calculated. This is the next result.
\begin{lema}\label{homint3} There is an isomorphism
$B_{\lambda^{-1}}\simeq \underline{\Hom}(\ku_{\lambda},
\ku_{\lambda})$ as $H$-module algebras.
\end{lema}
\pf First let us prove that there are natural isomorphisms
$$\Hom_H(X, B_{\lambda^{-1}})\simeq\, \Hom_A(X\ot_{\ku} \ku_{\lambda},\ku_{\lambda}) $$
for every $X\in \Rep(H)$. Let $\phi_X:\Hom_H(X,
B_{\lambda^{-1}})\to \Hom_A(X\ot_{\ku}
\ku_{\lambda},\ku_{\lambda}) $, and $\psi_X:$ $ \Hom_A(X\ot_{\ku}
\ku_{\lambda},\ku_{\lambda})$ $\to$$ \Hom_H(X, B_{-\lambda})$ be
the maps defined by
$$\phi_X(\alpha)(x\ot 1)=\alpha(x)(1), \quad \psi_X(\beta)(x)(t)= \beta(t\cdot x\ot 1), $$
for all $x\in X$, $t\in H$.

A straightforward computation shows that $\phi_X $ and $\psi_X$
are well defined maps one the inverse of the other. Now we must
prove that the algebra structure on the internal Hom described in
subsection \ref{subsection-modcat} coincides with the algebra
structure of $B_{\lambda^{-1}}$ given in \eqref{product1} via
these isomorphisms.

In this case is not hard to see that the evaluation map
$ev:B_{\lambda^{-1}}\ot_{\ku} \ku_{\lambda}\to \ku_{\lambda}$ is
$ev(f\ot 1)=\phi_{B_{\lambda^{-1}}}(\id)(f\ot 1)=f(1)$ for all
$f\in B_{\lambda^{-1}}$. The product
$\mu:B_{\lambda^{-1}}\ot_{\ku} B_{\lambda^{-1}}\to
B_{\lambda^{-1}}$ according to \ref{mult-intt} is
\begin{align*} \mu(f\ot g)(t)&= \psi_{B_{\lambda^{-1}}\ot_{\ku}
B_{\lambda^{-1}}}\big( ev\,(\id\ot ev)\,m_{B_{\lambda^{-1}},
B_{\lambda^{-1}},
\ku_{\lambda}}\big)(f\ot g)(t)\\
&=\big( ev\,(\id\ot ev)\,m_{B_{\lambda^{-1}}, B_{\lambda^{-1}},
\ku_{\lambda}}\big)((t\_1\cdot f\ot t\_2\cdot g)\ot 1)\\
&=\big( ev\,(\id\ot ev)\big)(J^{-1}(\lambda^{-1} )t\_1\cdot f\ot\, J^{-2}(\lambda^{-1})t\_2\cdot g\ot 1)\\
&=f(J^{-1}(\lambda^{-1})t\_1) \;
g(J^{-2}(\lambda^{-1})t\_2),\\
\end{align*}
for all $f, g\in B_{\lambda^{-1}}$, $t\in H$.
 \epf

\subsection{Dynamical Datum}\label{datum}
We introduce the notion of dynamical datum following ideas
contained in \cite{EN1}.

\begin{defi} A \emph{dynamical datum} for the pair $(H, A)$ is a
collection $(K, \{V_{\lambda}\}_{\lambda\in \widehat{A}})$ where
$K$ is an $H$-simple left $H$-comodule algebra, semisimple as an
algebra, such that $K^{\co H}=\ku$, and $V_{\lambda}$ is a
collection of irreducible $K$-modules, such that for every
$\lambda, \mu\in \widehat{A}$ there are $H$-module isomorphisms
\begin{equation}\label{dynamical-datum} \St_K(V_{\lambda},V_{\mu})
\simeq C(\lambda\,\mu^{-1}).
\end{equation}

\smallbreak

Two dynamical data $(K, \{V_{\lambda}\})$, $(S, \{W_{\lambda}\}) $
are \emph{equivalent} if and only if there exists an object
$P\in\, ^{H}\Mo_S$ such that $K\simeq \End_S(P_S)$ as $H$-comodule
algebras, and a family of $K$-module isomorphisms
$$\{\phi_{\lambda}:P\ot_S W_{\lambda}\to V_{\lambda}| \lambda\in
\widehat{A} \}.$$

\end{defi}

\begin{rmk} The $H$-module algebra structure on $\End_S(P_S)$ is
given in Lemma \ref{comodd}.
\end{rmk}

Our definition of dynamical datum is equivalent with the
definition given in \cite{EN1} when $H=\ku[G]$ is the group
algebra of a finite group $G$.

\medbreak

Indeed if $K$ is a $\ku[G]$-simple left comodule algebra such that
$K^{\co G}=\ku$ then $K=\ku^{\psi} [F]$ is the twisted group
algebra of a subgroup $F$ of $G$ and  $\psi\in Z^2(F,
\ku^{\times})$ is a normalized 2-cocycle.

If $(K, \{V_{\lambda}\})$ is dynamical datum in our sense then
$\St_{\ku^{\psi}[F]}(V_{\lambda},V_{\mu})\simeq
C(\lambda\,\mu^{-1})$. But we know that $\St_{\ku^{\psi}
F}(V_{\lambda},V_{\mu})\simeq \ku [G]\ot_{F}
\Hom(V_{\lambda},V_{\mu})$,  see \cite[Ex. 2.18]{AM}. Also,
$\Ind^G_A\, \lambda\simeq C(\lambda)$ as $\ku[G]$-modules, where
the isomorphism is given by $\beta:\Ind^G_A\, \lambda\to
C(\lambda)$, $\beta(\overline{g\ot
1})(h)=\delta_{g^{-1}}(P_{\lambda}\,h)$ for all $g, h\in G$.
Altogether implies that $ \ku G\ot_{\ku F} (V_{\mu}\ot_{\ku}
V^*_{\lambda}) \simeq \Ind^G_A\, (\lambda\,\mu^{-1})$, thus
$(\ku^{\hat{\psi}}[F], \{V^*_{\lambda}\})$ is a dynamical datum
according to the definition given in \cite{EN1}. Here
$\hat{\psi}(g,h)=\psi(h^{-1}, g^{-1})$ for all $g,h\in F$.

\medbreak

Is not hard to prove also that in this case two dynamical data
$(\ku^{\psi}[F], \{V_{\lambda}\})$, $(\ku^{\psi}[F'],
\{V'_{\lambda}\})$ are equivalent, according our definition, if
and only if there exists $g\in G$ such that $F'=ad_g F$ and the
corresponding representations are conjugated by $g$. This follows
from the fact that $\ku[G]$ is pointed and the quotient
$\ku^{G/F}$ is a cosemisimple coalgebra.

\begin{lema} Let $(K, \{V_{\lambda}\}_{\lambda\in \widehat{A}})$  be a dynamical datum.
Then\begin{enumerate}
    \item[1.]  $V_{\lambda}$ and $V_{\mu} $ are non-isomorphic $K$-modules
     if $\lambda\neq \mu$.
    \item[2.]  For all $\lambda\in\widehat{A}$
\begin{equation}\label{dim-dynamical} (\dim V_{\lambda})^2\,\mid A\mid= \dim
K.
\end{equation}
\end{enumerate}
\end{lema}
\pf (1) Assume that $\lambda\neq \mu$. Then
\begin{align*} \Hom_K(V_{\lambda}, V_{\mu})&\simeq \Hom_H(\ku,\St_K(V_{\lambda},
V_{\mu})) \simeq \Hom_H(\ku, C(\lambda\,\mu^{-1}))\\
&\simeq \Hom_H(\ku, (\Ind^H_{A} (\mu\,\lambda^{-1}))^*)\simeq
\Hom_{\ku A}(\mu\,\lambda^{-1},\ku)=0.
\end{align*}
The first isomorphism is a particular case of Proposition
\ref{stab-properties} (3) and the third is Remark \ref{isod}.

(2) By the definition of dynamical datum we have that
$\St_K(V_{\lambda})\simeq C(\varepsilon)$. Thus
$$\dim \St_K(V_{\lambda})=\frac{\dim H}{\mid A\mid}.$$
Now equation \eqref{dim-dynamical} follows from
\eqref{dim-stab}.\epf

\begin{rmk}\label{iso-clases} Observe that
\eqref{dim-dynamical} implies that the set
$\{V_{\lambda}\}_{\lambda\in \widehat{A}}$ is a complete set of
representatives of isomorphism classes of irreducible
representations of $K$.
\end{rmk}

\subsection{Dynamical twists constructed from a dynamical
datum}\label{exchange} Let $(K, \{V_{\lambda}\}_{\lambda\in
\widehat{A}})$ be a dynamical datum. We shall construct a
dynamical twist associated to $(K, \{V_{\lambda}\}_{\lambda\in
\widehat{A}})$. This procedure is called the \emph{exchange
construction} in \cite{EN1}, see also \cite{EV}.

\medbreak

If $X\in \Rep(H)$ and $V\in {}_K\Mo$, we always assume that the
vector space $X\ot_{\ku}  V$ carries the left $K$-action described
in \eqref{mod-action}.

\medbreak

For any pair $\lambda, \mu\in \widehat{A}$ choose $H$-module
isomorphisms
$$\omega_{\lambda, \mu}:\St_K(V_{\lambda},V_{\mu})\simeq
C(\lambda\,\mu^{-1})$$ such that for $\lambda=\mu$ the identity of
$\St_K(V_{\lambda})$ is mapped to $\varepsilon$.

\medbreak

For any $\lambda,\mu \in \widehat{A}$, $x\in X[\mu]$ denote by
$\Psi(\lambda, x):X^*\ot_{\ku} V_{\lambda}\to V_{\lambda\,\mu}$
the $K$-map obtained as the image of $x$ under the (natural)
isomorphisms
\begin{equation}\label{chain-iso}
\begin{split} X[\mu]&\simeq\Hom_{A}(\mu, \Res^H_A X)
\simeq \Hom_{H}(\Ind^H_{A}\, \mu, X)\\ &\simeq\Hom_{H}(X^*,
(\Ind^H_{A}\, \mu)^*) \simeq
\Hom_{H}(X^*,C(\mu^{-1}))\\&\simeq\Hom_{H}(X^*,\St_K(V_{\lambda},V_{\lambda\,\mu}))\simeq
 \Hom_{K}(X^*\ot_{\ku} V_{\lambda},V_{\lambda\,\mu}).
\end{split}
\end{equation}
The second isomorphism by Frobenius reciprocity \ref{frobenius},
and the fourth is remark \ref{isod}, the fifth isomorphism is
defined by $\omega_{\lambda,\lambda\,\mu}$ and the last
isomorphism comes from Proposition \ref{stab-properties} (3).

\medbreak

More explicitly if $f\in X^*$, $v\in V_{\lambda}$ then
\begin{align}\label{def-psi} \Psi(\lambda, x)(f\ot v)=
\omega^{-1}_{\lambda,\,\lambda\,\mu}(\widetilde{f}^x)(1\ot v),
\end{align}
where $\widetilde{f}^x$ is the element in $C(\mu^{-1})$ defined by
$\widetilde{f}^x(h)=f(\Ss(h)\cdot x),$ for all $h\in H$.

\begin{lema} Let $X, Y\in \Rep(H)$ and $f:X\to Y$ an
$H$-module map. If $\lambda, \mu\in \widehat{A} $ and $x\in
X[\mu]$ then
\begin{align}\label{psi} \Psi(\lambda, f(x))=\Psi(\lambda,
x)(f^*\ot\,\id_{V_{\lambda}}).
\end{align}
\end{lema}

\pf Straightforward. It follows from the naturality of the
isomorphisms \eqref{chain-iso} or directly from \ref{def-psi}.
\epf

For any $X, Y\in \Rep(H)$, $V\in {}_K\Mo$ we shall denote by
$\phi_{XY}:( X\ot_{\ku} Y)^*\to Y^*\ot_{\ku} X^*$ the canonical
isomorphism and $m_{XYV}:(X\ot_{\ku} Y)\ot_{\ku} V\to X\ot_{\ku}
(Y\ot_{\ku} V)$ the canonical associativity isomorphism.

\medbreak

Let $Y$ be another $H$-module and $y\in Y([\eta]$. The composition
\begin{align*}&( X\ot_{\ku} Y)^*\ot_{\ku} V_{\lambda}\xrightarrow{\phi_{XY}\ot\id}(Y^*\ot_{\ku} X^*)\ot_{\ku}
V_{\lambda}\xrightarrow{\,\,m_{Y^*,X^*,V_{\lambda}}\,\,} {Y^*\ot_{\ku} (X^*\ot_{\ku} V_{\lambda})}\longrightarrow \\
&\xrightarrow{\,\,\id\ot \Psi(\lambda,x)\,\,}{Y^*\ot_{\ku}
V_{\lambda\,\mu}} \xrightarrow{\,\,\Psi(\lambda\,\mu, y)\,\,}{
V_{\lambda\,\mu\,\eta}}
\end{align*}
determines a unique element in $(X\ot_{\ku} Y)[\mu\eta]$, that we
denote by $I_{XY}(\lambda)(x\ot y)$. That is, we have defined a
map $I_{XY}(\lambda):X\ot_{\ku} Y\to X\ot_{\ku} Y$ by
\begin{align}\label{twist-dynam1} \Psi(\lambda,I_{XY}(\lambda)(x\ot
y))= \Psi(\lambda\,\mu, y)\, (\id\ot
\Psi(\lambda,x))\,m_{Y^*,X^*,V_{\lambda}}(\phi_{XY}\ot\id).
\end{align}

By \eqref{psi} the maps $I_{XY}$ are natural, in particular, there
is an element $I(\lambda)\in H\ot_{\ku} H$ such that
\begin{align}\label{naturalidad} I_{XY}(\lambda)(x\ot
y)=I(\lambda)(x\ot y),
\end{align}
for all $\lambda\in \widehat{A}$, $x\in X$, $y\in Y$.

\begin{lema}\label{dynamic21}
$I(\lambda)\in H\ot_{\ku} H$ is an invertible element for all
$\lambda\in \widehat{A}$.
\end{lema}

\pf The proof is entirely similar to the proof of \cite[Lemma
6.2]{EN1}. For completeness we will write it down.  By the
definition \eqref{twist-dynam1} of $I_{XY}$ the surjectivity of
the map
$$I_{XY}(\lambda):\bigoplus_{\nu} X[\mu\nu^{-1}]\ot_{\ku} Y[\nu\eta]\to
(X\ot_{\ku} Y)[\mu\eta]$$ is equivalent to the surjectivity of the
composition
$$\bigoplus_{\nu}\Hom_K(V_{\nu}, Y\ot_{\ku} V_{\eta}) \ot_{\ku}
\Hom_K(X^*\ot_{\ku} V_{\lambda}, V_{\nu})\to \Hom_K(X^*\ot_{\ku}
V_{\lambda},Y\ot_{\ku} V_{\eta})
$$
and this follows since $X^*\ot_{\ku} V_{\lambda}\simeq
\bigoplus_{\nu} X[\nu\lambda^{-1}]\ot_{\ku} V_{\nu}$; indeed by
the isomorphisms \eqref{chain-iso} $\Hom_K(X^*\ot_{\ku}
V_{\lambda}, V_{\nu})$ $\simeq X[\nu\lambda^{-1}],$  so a copy of
$\bigoplus_{\nu} X[\nu\lambda^{-1}]\ot_{\ku} V_{\nu}$ is inside of
$X^*\ot_{\ku} V_{\lambda}$ but
\begin{align*}\dim (\bigoplus_{\nu} X[\nu\lambda^{-1}]\ot_{\ku} V_{\nu})&=\sum_{\nu}
\dim(X[\nu\lambda^{-1}]) \dim(V_{\nu})\\
&=\sum_{\nu} \dim(X[\nu\lambda^{-1}]) \dim(V_{\lambda})\\
&= \dim(X)\dim(V_{\lambda})= \dim(X^*\ot_{\ku} V_{\lambda}).
\end{align*}
Here the second equality follows from \eqref{dim-dynamical}, that
is all  modules $V_{\nu}$ have the same dimension.

 \epf

Let us define $J:\widehat{A}\to H\ot_{\ku} H$ by
$J(\lambda)=\Ss^{-1}(I^{-2}(\lambda))\ot\,\Ss^{-1}(I^{-1}(\lambda))$.
\begin{prop}\label{dynamic2} $J(\lambda)$ is a dynamical twist for
$H$.

\end{prop}
\pf Let $X, Y, Z\in\Rep(H)$, $x\in X(\mu)$, $y\in Y(\eta)$, $z\in
Z(\nu)$ and $\mu, \eta,\nu, \lambda\in \widehat{A}$. If $a\in A$
then
$$I_{XY}(\lambda)(a\cdot x\ot \,a\cdot y)=(\mu\eta)(a)I_{XY}(\lambda)(x\ot y)
=a\cdot  I_{XY}(\lambda)(x\ot y).$$ Hence $I(\lambda)$ and
therefore $J(\lambda)$ commutes with $a\ot\, a$ for all $a\in A$.

\medbreak

Set $m_1= m_{Z^*,(X\ot Y)^*,V_{\lambda} }$, $m_2=m_{(Y\ot Z)^*,
X^*, V_{\lambda}}$. The associativity implies that
$$\Psi(\lambda\mu\eta,z)\,\big(\id_{Z^*}\ot
\Psi(\lambda, I_{XY}(\lambda) (x\ot y))\big) m_1 (\phi_{X\ot Y,
Z}\ot\id)$$ equals to
$$\Psi(\lambda\mu, I_{YZ}(\lambda\mu)(y\ot z))\,\left(\id_{(Y\ot Z)^*}\ot
\Psi(\lambda,x)  \right)m_2(\phi_{X,Y\ot Z}\ot\id)(a^*\ot\id),$$
thus $I_{X\ot Y,Z}(\lambda)( I_{XY}(\lambda)(x\ot y)\ot z)= I_{X,
Y\ot Z}(\lambda)(x\ot I_{YZ}(\lambda\mu)(y\ot z))$. This implies
that
$$(\Delta\ot\id)I(\lambda)\, (I(\lambda)\ot 1)= \sum_{\mu}(\id\ot\Delta)I(\lambda)
\,(P_{\mu}\ot I(\lambda\mu)).  $$ Using the definition of $J$, the
fact that $\Ss(P_{\mu})=P_{\mu^{-1}}$ and properties of the
antipode we get that $J(\lambda)$ satisfies equation \ref{d2}.
Also, since we have chosen isomorphisms $w_{\lambda,\lambda}$ such
that maps the identity of $\St_K(V_{\lambda})$ to $\varepsilon$
then $J(\lambda)$ verifies identity \eqref{d3}.

\epf

We shall say that $J:\widehat{A}\to H\ot H$ is the dynamical twist
\emph{associated} to the dynamical datum $(K, \{V_{\lambda}\})$.

\begin{prop}\label{mod-datum} Let $(K, \{V_{\lambda}\})$ be a
dynamical datum and $J:\widehat{A}\to H\ot H$ the dynamical twists
associated. Then
$$ {}_K\Mo \simeq \Mo^{(J)}$$
as module categories over $\Rep(H)$.
\end{prop}
\pf We know that for any $\lambda\in \widehat{A}$ there are module
equivalences $ {}_K\Mo\simeq \Rep(H)_{\St_K(V_{\lambda})}$ and
$\Mo^{(J)}\simeq \Rep(H)_{B_{\lambda}}$. The first one is
Corollary \ref{modk} and the second follows from Lemma
\ref{homint3} and \cite[Theorem 3.17]{eo}. We shall prove that the
$H$-module isomorphism
$w_{\lambda}=w_{\lambda,\lambda}:\St_K(V_{\lambda})\to
B_{\lambda}$ is an algebra isomorphism. This will end the proof.

\medbreak

Let $\lambda\in \widehat{A},$ $v\in V_{\lambda}$ and $f, g\in
B_{\lambda}$. Denote $X=H\diagup (\ku A)^+ H$, thus $X^*=
C(\varepsilon)$. Then using \eqref{def-psi} we get
\begin{align*} \Psi(\lambda, \overline{1})\,
(\id\ot \Psi(\lambda,\overline{1}))(f\ot g\ot v)&=\Psi(\lambda,
\overline{1})(g\ot\,
\omega^{-1}_{\lambda}(\widetilde{f})(1\ot v))\\
&=\omega^{-1}_{\lambda}(\widetilde{g})\;\omega^{-1}_{\lambda}(\widetilde{f})(1\ot
v)\\
&=\omega^{-1}_{\lambda}(g\circ \Ss)\;\omega^{-1}_{\lambda}(f\circ
\Ss)(1\ot v)
\end{align*}
Recall that here $\widetilde{f} \in B_{\lambda}$ denotes the map
$\widetilde{f}(h)=f(\Ss(h))$, for all $h\in H$. On the other hand
\begin{align*}  \Psi(\lambda,I_{XX}(\lambda)(\overline{1}\ot
\overline{1}))(f\ot g\ot
v)=\omega^{-1}_{\lambda}(\widetilde{f}^{I^1(\lambda)}\ot
\widetilde{g}^{I^2(\lambda)})(1\ot v),
\end{align*}
where, for all $h\in H$
\begin{align*}(\widetilde{f}^{I^1(\lambda)}\ot\,
\widetilde{g}^{I^2(\lambda)})(h)&= (f\ot g)(\Ss(h)\cdot I(\lambda))\\
&=(g\circ \Ss)(J^{-1}(\lambda)h\_1)(f\circ
\Ss)(J^{-2}(\lambda)h\_2)\\ &=(g\circ \Ss)\cdot (f\circ \Ss)(h).
\end{align*}
The product on last equality is the product in $B_{\lambda}$.
Therefore $$\omega^{-1}_{\lambda}(g\circ
\Ss)\;\omega^{-1}_{\lambda}(f\circ
\Ss)=\omega^{-1}_{\lambda}((g\circ \Ss)\cdot (f\circ \Ss)),$$ and
this ends the proof.
 \epf

The construction of the dynamical twist from the dynamical datum
is not canonical, however, in the following we shall prove that
equivalent dynamical data defines the same gauge equivalence class
of dynamical twist. First we need the next technical Lemma.

\medbreak

Let $K, S$ be left $H$-comodule algebras and $P\in \,
^{H\!}_{K\!}\Mo_S$. For any $X\in \Rep(H)$, $M\in {}_S\Mo$ let
$$\theta_{X,M}:X\ot_{\ku}( P\ot_S M)\to P\ot_S(X \ot_{\ku} M)$$ be
defined by $\theta_{X,M}(x\ot p\ot\, m)=p\_0\ot\,
\Ss^{-1}(p\_{-1})\cdot x\ot\, m, $ for any $x\in X, m\in M, p\in
P$.
\begin{lema} Let $X, Y\in \Rep(H),$ $M, N \in {}_S\Mo$, then
the maps $\theta_{X,M}$ are well-defined $K$-isomorphisms. Also if
$g:M\to N$ is an $S$-module map and $f:X\to Y$ an $H$-module map
we have
\begin{align}\label{theta1} \theta_{X\ot_{\ku} Y,M}&=\theta_{X,Y\ot_{\ku}
M}(\id_X\ot\, \theta_{X,M}),\\
\label{theta2} (\id_P\ot\id_X\ot g)\, \theta_{X, M}&=
\theta_{X,N}\, (\id_X\ot\id_P\ot g),\\
\label{theta3} (\id_P\ot f\ot \id_M) \, \theta_{X,M}&=
\theta_{Y,M} (f\ot \id_P\ot \id_M).
\end{align}
\end{lema}
\pf Straightforward.

\epf

\begin{prop}\label{gauge-t} If $(K, \{V_{\lambda}\})$ and $(S, \{W_{\lambda}\})$ are
equivalent dynamical data then the associated dynamical twists are
gauge equivalent.
\end{prop}
\pf Let $J, J'$ be the dynamical twists associated to $(K,
\{V_{\lambda}\})$ and $(S, \{W_{\lambda}\})$ respectively and
correspondingly the maps $\Psi, \Psi'$ and $I, I'$. Since $(K,
\{V_{\lambda}\})$ and $(S, \{W_{\lambda}\})$ are equivalent there
exists $P\in\, ^{H}\Mo_S$ such that $K\simeq \End_S(P_S)$ as
$H$-comodule algebras, and $K$-module isomorphisms
$\phi_{\lambda}:P\ot_S W_{\lambda}\to V_{\lambda}$.

Let  $\mu, \lambda\in \widehat{A}$ and $X\in \Rep(H)$. For any
$x\in X[\mu]$,  define $\sigma_X(\lambda)(x)$ the element obtained
as the preimage of the map

\begin{align*} &X^*\ot_{\ku} V_{\lambda}\xrightarrow{\;\id_{X^*}\ot \phi^{-1}_{\lambda}\;}
 X^*\ot_{\ku} (P\ot_S W_{\lambda})\xrightarrow{\;\theta_{X^*,W_{\lambda}}\;}
 P\ot_S (X^*\ot_{\ku} W_{\lambda})\longrightarrow\\
&\xrightarrow{\;\id_P\ot \Psi'(\lambda,x)\;} P\ot_S
W_{\lambda\mu}\xrightarrow{\;\phi_{\lambda\mu}\;} V_{\lambda\mu}.
\end{align*}
under the ismorphisms \ref{chain-iso}. That is
\begin{align}\label{gauge-data} \Psi(\lambda, \sigma_X(\lambda)(x))=
\phi_{\lambda\mu}\, (\id_P\ot
\Psi'(\lambda,x))\,\theta_{X^*,W_{\lambda}}\,(\id_{X^*}\ot\,\phi^{-1}_{\lambda}).
\end{align}

We claim that the maps $\sigma_X(\lambda)$ are natural
isomorphisms. Indeed, if  $x\in X[\mu]$, $Y\in \Rep(H)$ and
$f:X\to Y$ is an $H$-module map then
\begin{align*}\Psi(\lambda, (\sigma_Y(\lambda)\circ f)(x))&=
\Psi(\lambda, \sigma_Y(\lambda)(f(x)))\\
&=\phi_{\lambda\mu}\, (\id_P\ot \Psi'(\lambda,f(x)))\,\theta_{Y^*,W_{\lambda}}\,(\id_Y^*\ot\,\phi^{-1}_{\lambda})\\
&=\phi_{\lambda\mu}\, \big(\id_P\ot
\Psi'(\lambda,x)(f^*\ot\id_{W_{\lambda}})
\big)\,\theta_{Y^*,W_{\lambda}}\,(\id_Y^*\ot\,\phi^{-1}_{\lambda})\\
&=\phi_{\lambda\mu}\, \big(\id_P\ot
\Psi'(\lambda,x)\big)\theta_{X^*,W_{\lambda}} (f^*\ot\,
\phi^{-1}_{\lambda})\\
&=\Psi(\lambda, \sigma_X(\lambda)(x)) (f^*\ot\, \id_{V_{\lambda}})\\
&=\Psi(\lambda, (f\circ \sigma_X(\lambda))(x)).
\end{align*}
The third equality by \eqref{psi}, the fourth by \eqref{theta3}
and the sixth again by \eqref{psi}. Therefore there exists an
invertible element $\sigma(\lambda)\in H$  such that
$$\sigma_X(\lambda)(x)=\sigma(\lambda)\cdot x$$ for all $\lambda \in
\widehat{A}, x\in X$. Set $t(\lambda)=
\Ss^{-1}(\sigma(\lambda)^{-1})$. We shall prove that $t(\lambda)$
is a gauge equivalence between $J$ and $J'$. Clearly $t(\lambda)$
verifies \eqref{gauge3} and \eqref{gauge1}. Let us prove that
\eqref{gauge2} holds. Let $\eta\in \widehat{A}$, $Y\in \Rep(H)$,
$y\in Y[\eta]$ then
\begin{align*}&\Psi(\lambda, \sigma_{X\ot
Y}(\lambda)I'_{XY}(\lambda)(x\ot y)))=\\&= \phi_{\lambda\mu\eta}\,
(\id_P\ot \Psi'(\lambda,I'_{XY}(x\ot
y)))\,\theta_{(X\ot Y)^*,W_{\lambda}}\, (\id_{(X\ot Y)^*}\ot\, \phi^{-1}_{\lambda})\\
&=\phi_{\lambda\mu\eta}\,(\id_P\ot \Psi'(\lambda\mu, y))(\id_{P\ot
Y^*}\ot \Psi'(\lambda, x))\,\theta_{(X\ot Y)^*,
W_{\lambda}}\,(\id_{(X\ot Y)^*}\ot\,\phi^{-1}_{\lambda})\\
&=\Psi(\lambda, \sigma_Y(\lambda\mu)(y))\big(\id_{Y^*}\ot\,
\Psi'(\lambda, \sigma_X(\lambda)(x)\big)(\id_{(X\ot Y)^*}\ot\,
\phi_{\lambda}))\\
&(\id_{Y^*}\ot \,\theta^{-1}_{X^*,W_{\lambda}})\,
\theta^{-1}_{Y^*,W_{\lambda}}\,\theta_{(X\ot Y)^*,
W_{\lambda}}\,(\id_{(X\ot Y)^*}\ot\,\phi^{-1}_{\lambda})\\
&=\Psi(\lambda, I_{XY}(\lambda)(\sigma_X(\lambda)(x)\ot\,
\sigma_Y(\lambda\mu)(y)))
\end{align*}

The third equality follows from \eqref{theta2} since
$\Psi'(\lambda, x)$ is an $S$-module morphism and the fourth
follows from  \eqref{theta1}. Thus
$$ I_{XY}(\lambda)(\sigma_X(\lambda)(x)\ot\,
\sigma_Y(\lambda\mu)(y))= \sigma_{X\ot
Y}(\lambda)I'_{XY}(\lambda)(x\ot y),$$ and this implies
\eqref{gauge2}. \epf

\subsection{Dynamical datum constructed from a dynamical
twist}\label{dyn} \

Let $J:\widehat{A}\to H \ot H$ be a dynamical twist and consider
the module category $\Mo^{(J)}$ explained in subsection
\ref{modcat-twist}. By Theorem \ref{mod-overhopf} there exists an
$H$-simple left $H$-comodule algebra $K$, with $K^{\co H}=\ku$ and
a module equivalence $(F,c):\Mo^{(J)}\to {}_K \Mo$. Set
$V_{\lambda}=F(\ku_{\lambda^{-1}})$ for any
$\lambda\in\widehat{A}$.

\begin{prop}\label{dynamdat} The pair $(K,\{V_{\lambda}\})$ is a
dynamical datum.
\end{prop}

\pf Let $X\in \Rep(H)$, $\lambda, \mu\in \widehat{A}$. Then
\begin{align*} &\Hom_H(X, \St_K(V_{\lambda},V_{\mu}))\simeq
\Hom_K(X\ot_{\ku} V_{\lambda}, V_{\mu})\simeq \\&\simeq
\Hom_K(X\ot_{\ku} F(\ku_{\lambda^{-1}}), F(\ku_{\mu^{-1}}))\simeq
\Hom_K(F(X\ot_{\ku}\, \ku_{\lambda^{-1}}),
F(\ku_{\mu^{-1}}))\simeq\\&\simeq \Hom_A(X\ot_{\ku}\,
\ku_{\lambda^{-1}}, \ku_{\mu^{-1}})\simeq \Hom_A(\Res^H_A
X,\ku_{\lambda\mu^{-1}})\simeq\\
&\simeq \Hom_A(\ku_{\mu\lambda^{-1}},\Res^H_A X^*)\simeq
\Hom_H(\Ind^H_A (\mu\lambda^{-1}), X^*) \simeq\\
&\simeq  \Hom_H( X,(\Ind^H_A (\mu\lambda^{-1})^*)\simeq \Hom_H(
X,C(\lambda\mu^{-1})).
\end{align*}
The last isomorphism by Remark \ref{isod}. Thus by Yoneda's Lemma
there is an $H$-module isomorphism $
\St_K(V_{\lambda},V_{\mu})\simeq C(\lambda\mu^{-1})$ . \epf

The equivalence class of the dynamical data constructed from a
dynamical twist as above does not depend on the gauge equivalence
class of the dynamical twists.  This is evident from Proposition
\ref{eq1}.

\subsection{Main result}\label{main0} Using the same notation as in \cite{EN1}
we shall denote by $T$ and $D$ the maps between gauge equivalence
classes of dynamical twists and equivalence classes of dynamical
data described in subsections \ref{exchange} and \ref{dyn}
respectively. That is,

\begin{equation*}
\left\{
\begin{array}{c}
\mbox{gauge equivalence}\\
\mbox{classes of} \\
\mbox{dynamical twists} \\
J: \widehat{A}\to H\otimes H
\end{array}
\right\}
\begin{array}{c}
{}\\
\xrightarrow{\quad\; D \quad \;}\\
\xleftarrow{\quad\; T \quad\;}\\
{}
\end{array}
\left\{
\begin{array}{c}
\mbox{equivalence}\\
\mbox{classes of} \\
\mbox{dynamical data}\\
 (K,\,\{ V_\lambda\}_{\lambda\in \widehat{A}} )
\end{array}
\right\}.
\end{equation*}

\begin{teo}\label{main} The maps $D$ and $T$ are inverses of each
other.
\end{teo}
\pf First we shall prove that $D\circ T= \Id$. Let
$(K,\{V_{\lambda}\})$ be a dynamical data and $J: \widehat{A}\to
H\ot H$ the dynamical twist coming from the exchange construction
according to subsection \ref{exchange}. By Proposition
\ref{mod-datum} ${}_K\Mo\simeq \Mo^{(J)}$. Let
$(S,\{W_{\lambda}\})$ be a dynamical data as constructed in
subsection \ref{dyn}. By definition ${}_S\Mo\simeq \Mo^{(J)}$,
then ${}_S\Mo\simeq {}_K\Mo$, therefore, using Proposition
\ref{eq1}, $(S,\{W_{\lambda}\})$ is equivalent to
$(K,\{V_{\lambda}\})$.

\medbreak

Now, let us prove that $T\circ D= \Id$. Let $J_1: \widehat{A}\to
H\ot H$ be a dynamical twist and $(K,\{V_{\lambda}\})$ be the
dynamical data constructed as in subsection \ref{dyn}. In
particular this means that there is a module equivalence
$(F,c):\Mo^{(J_1)}\to {}_K\Mo$ such that
$F(\ku_{\lambda^{-1}})=V_{\lambda}$ for all $\lambda\in
\widehat{A}$. Let $J_2: \widehat{A}\to H\ot H$ be the dynamical
data associated to $(K,\{V_{\lambda}\})$. Let $I_1,
I_2:\widehat{A}\to H\ot H$ be the maps defined as
$$ I_i(\lambda)=\Ss(J_i^{-2}(\lambda))\ot \,\Ss(J_i^{-1}(\lambda)), \,\, i=1,2.$$
Let $\Psi_2$ be the map defined form the dynamical data
$(K,\{V_{\lambda}\})$ as in \eqref{def-psi}. By the exchange
construction we know that
\begin{align}\label{t2} \Psi_2(\lambda,I_2(\lambda)\cdot (x\ot
y))= \Psi_2(\lambda\mu, y)\, (\id\ot
\Psi_2(\lambda,x))\,m_{Y^*,X^*,V_{\lambda}}(\phi_{XY}\ot\id),
\end{align}
for all $\lambda,\mu,\eta\in\widehat{A}$, $X, Y\in \Rep(H)$, $x\in
X[\mu]$, $y\in Y[\eta]$.

For any $ \mu\in\widehat{A}$, $X\in \Rep(H)$, $x\in X[\mu]$ we
denote by $\Psi_1(\lambda,x):X^*\ot \ku_{\lambda^{-1}}\to
\ku_{(\mu\lambda)^{-1}}$ the map obtained as the image of $x$
under the composition of isomorphisms
$$X[\mu]\simeq \Hom_A(\ku_{\mu}, X)\simeq
\Hom_A(X^*, \ku_{\mu^{-1}})\simeq\Hom_A(X^*\ot_{\ku}\,
\ku_{\lambda^{-1}}, \ku_{(\mu\lambda)^{-1}}). $$

An easy computation shows that for all $x\in X[\mu]$,  $y\in
Y[\eta]$
$$\Psi_1(\lambda,I_1(\lambda)\cdot (x\ot
y))= \Psi_1(\lambda\mu, y)\big(\id\ot\Psi_1(\lambda,x)\big)
\widetilde{m}_{Y^*,X^*,\ku_{-\lambda}}(\phi_{XY}\ot\id),
$$
where $\widetilde{m}$ is the associativity as in \eqref{assoc-mc}.
\medbreak

To prove that $J_1$ is gauge equivalent to $J_2$ we will use the
same idea as in Proposition \ref{gauge-t}. For any $x\in X[\mu]$
define the maps $\sigma_X(\lambda):X\to X$ by
\begin{align*}\Psi_2(\lambda,\sigma_X(\lambda)(x))=F(\Psi_1(\lambda,x))\;
c^{-1}_{X^*,\ku_{\lambda^{-1}}}.
\end{align*}
The maps $\sigma_X(\lambda)$ are natural isomorphisms, hence there
exists an invertible element $\sigma(\lambda)\in H$ such that
$\sigma_X(\lambda)(x)=\sigma(\lambda)\cdot x$. Set $t(\lambda)=
\Ss^{-1}(\sigma(\lambda)^{-1})$. We  will prove that $t(\lambda)$
defines a gauge equivalence between $J_1$ and $J_2$. Now
\begin{align*}&\Psi_2(\lambda,\sigma_{X\ot Y}(\lambda)I_1(\lambda)\cdot (x\ot y))
= F(\Psi_1(\lambda,I_1(\lambda)\cdot (x\ot y)))\,c^{-1}_{(X\ot
Y)^*,\ku_{\lambda^{-1}}}\\
&=F(\Psi_1(\lambda\mu,
y))F(\id\ot\Psi_1(\lambda,x))F(\widetilde{m}_{Y^*,X^*,\ku_{\lambda^{-1}}})
\,c^{-1}_{Y^*\ot
X^*,\ku_{\lambda^{-1}}}\\
&=F(\Psi_1(\lambda\mu, y))F(\id\ot\Psi_1(\lambda,x))\,
c^{-1}_{Y^*, X^*\ot\ku_{\lambda^{-1}}}(\id\ot
c^{-1}_{X^*,\ku_{\lambda^{-1}}})\, m_{Y^*,X^*,V_{\lambda}}\\
&=F(\Psi_1(\lambda\mu^{-1},
y))\,c^{-1}_{Y^*,\ku_{(\mu\lambda)^{-1}}}(\id\ot
F(\Psi_1(\lambda,x))\, c^{-1}_{X^*,\ku_{\lambda^{-1}}})\,
m_{Y^*,X^*,V_{\lambda}}\\
&=\Psi_2(\lambda\mu, \sigma_Y(\lambda\mu)(y)) (\id\ot
\Psi_2(\lambda, \sigma_X(\lambda)(x)))\, m_{Y^*,X^*,V_{\lambda}}\\
&=\Psi_2(\lambda, I_2(\lambda)\cdot (\sigma_X(\lambda)(x))\ot\,
\sigma_Y(\lambda\mu)(y))).
\end{align*}
The third equality by \eqref{modfun1} and the fourth by the
naturality of $c$. Therefore
$$I_2(\lambda)\cdot (\sigma_X(\lambda)(x))\ot\,
\sigma_Y(\lambda\mu)(y)=\sigma_{X\ot Y}(\lambda)I_1(\lambda)\cdot
(x\ot y),$$ and this equality implies that $t(\lambda)$ satisfies
\eqref{gauge2}.
 \epf

\begin{rmk}\label{usual-twist} As a immediate consequence of Theorem \ref{main} we note that
gauge equivalence classes of (usual) twists for Hopf algebras are
parameterized by equivalence classes of pairs $(K,V)$ where
\begin{itemize}
    \item $K$ is a semisimple $H$-simple left $H$-comodule
    algebra,
    \item $K^{\co H}=\ku$, and
    \item $\St_K(V)\simeq H^*$ as left
$H$-modules
\end{itemize}
If this is the case then $K$ is simple algebra, something expected
since the category ${}_K\Mo$ must have only one simple object.
\end{rmk}

\section{Some examples }\label{examples}

In this section we shall give some examples of dynamical data and
we compute the corresponding dynamical twist.


\subsection{Case when $K=\ku[A]$}\label{k=a}
\

This example is \cite[Ex. 6.10]{EN1} for an arbitrary Hopf
algebra. Let $K=\ku[A]$ the group algebra of the group $A$. Let
$f:\widehat{A}\to \widehat{A}$ be a bijection. Assume that for all
$\lambda, \mu \in \widehat{A}$ there exists an element
$g(\lambda,\mu)$ in the normalizer $N(A)$ of $A$ such that
\begin{align}\label{normal}
\big(f(\lambda)f(\mu)^{-1}\big)(a)=(\mu\lambda^{-1})(g(\lambda,\mu)\,a
\, g(\lambda,\mu)^{-1})
\end{align}
for all $a\in A$. Set $V_{\lambda}=\ku_{f(\lambda)}$, then $(K,
\{V_{\lambda}\})$ is a dynamical datum for $(H,A)$.
\pf Clearly $K$ is an $H$-simple left $H$-comodule algebra,
semisimple as an algebra with trivial coinvariants. By Proposition
\ref{hopf-galois} $\St_{\ku[A]}(V_{\lambda}, V_{\mu})\simeq
\Hom_{\ku[A]}(H, \Hom_{\ku}(V_{\lambda}, V_{\mu})$. Define the map
$\omega_{\lambda,\mu}:\Hom_{\ku[A]}(H,
\ku_{f(\mu)f(\lambda)^{-1}})\to C(\lambda\mu^{-1})$ as follows. If
$\alpha\in \Hom_{\ku[A]}(H, \ku_{f(\mu)f(\lambda)^{-1}})$, $h\in
H$ then
$$\omega_{\lambda,\mu}(\alpha)(h)=\alpha(g(\lambda,\mu)^{-1}h).$$
Equation \eqref{normal} implies that
$\omega_{\lambda,\mu}(\alpha)\in  C(\lambda\mu^{-1})$.\epf

Now we compute the corresponding dynamical twist associated to
$(K, \{V_{\lambda}\})$. If $X\in \Rep(H), x\in X[\mu]$, in this
case the maps $\Psi:X^*\ot_{\ku}\, \ku_{f(\lambda)}\to
\ku_{f(\lambda\mu)}$ are
$$ \Psi(\lambda,x)(f\ot 1)=f(g(\lambda,\lambda\mu)^{-1}\cdot x)$$
for all $f\in X^*$. Let us compute the corresponding dynamical
twist. Let $\mu, \eta \in \widehat{A}$ and $f_1, f_2 \in H^*$,
then $\Psi(\lambda, I(\lambda)(P_{\mu}\ot\, P_{\eta}))(f_1\ot
f_2\ot 1)$ is equal to
\begin{align*}
f_1(g(\lambda,\lambda\mu\eta)^{-1}I^2(\lambda)
P_{\eta})\;f_2(g(\lambda,\lambda\mu\eta)^{-1}I^1(\lambda)P_{\mu}).
\end{align*}
On the other hand we have that $$\Psi(\lambda\mu,
P_{\eta})(\id_{H^*}\ot \Psi(\lambda, P_{\mu}) (f_1\ot f_2\ot 1)$$
is equal to
$$f_1(g(\lambda,\lambda\mu)^{-1} P_{\mu})\;
 f_2(g(\lambda,\lambda\eta)^{-1} P_{\eta}).$$
Hence $$I^{-1}(\lambda)=\sum_{\mu, \eta} g(\lambda,\lambda\mu)\,
g(\lambda,\lambda\mu\eta)^{-1} P_{\mu}\ot\,
g(\lambda,\lambda\eta)\,g(\lambda,\lambda\mu\eta)^{-1} P_{\eta}.$$
Thus the dynamical twist in this case is
$$J(\lambda)=\sum_{\mu, \eta} P_{\mu}g(\lambda,\lambda\mu^{-1}\eta)
g(\lambda,\lambda\mu^{-1})^{-1}\ot
P_{\eta}g(\lambda,\lambda\mu^{-1}\eta)
g(\lambda,\lambda\eta^{-1})^{-1}.$$

\subsection{Dynamical twists for the Taft Hopf
algebras}\label{taft-ex}

In this subsection for each $c\in \ku^{\times}$ we construct a
dynamical twist for the Taft Hopf algebras.

\medbreak

Let $q$ be a $n$-primitive root of 1. Recall that the Taft algebra
$T(q)$ is the algebra generate by $g, x$ subject to the relations
$x^n=0, g^n=1$, $gx= q\, xg$. The Hopf algebra structure is
determined by
$$ \Delta(g)=g\ot g,\; \Delta(x)= 1\ot x+ x\ot g,\; \varepsilon(g)=1,\;
\varepsilon(x)=0,$$
$$\Ss(g)=g^{-1},\; \Ss(x)= -xg^{-1}. $$

Let $d\in \Na$ a divisor of $n$ and set $n=dm$. For any $c\in \ku$
denote by $\Ac(d,c)$ the algebra generated by $h$ and $y$ subject
to the relations $y^n=c.1$, $h^d=1$ and $hy=q^m\, yh$.

Define $\delta:\Ac(d,c)\to T(q)\ot_{\ku}\, \Ac(d,c)$ by
\begin{equation*} \delta(h)= g^{-m}\ot h, \;\,\;\,
\delta(y)=g^{-1}\ot\, y - xg^{-1}\ot 1.
\end{equation*}

These algebras are left $T(q)$-module algebra $T(q)$-simple.
Moreover $\Ac(d,c)$ is semisimple if and only if $c\neq 0$. This
algebras were considered in \cite{MoSch}, and also in \cite{eo}
where they classify indecomposable exact module categories over
$\Rep(T(q))$. \medbreak

Fix $c, b\in \ku^{\times}$ such that $b^n=c$. For $i=0\dots n-1$
let $V_i$ be the one-dimensional vector space generated by $v_i$
together with an action of $\Ac(1,c)$ defined by $ y\cdot
v_i=q^ib\, v_i$. The collection $\{V_i\}$ is a complete set of
representatives of isomorphism classes of irreducible modules of
$\Ac(1,c)$.

\medbreak

We shall prove that $(\Ac(1,c), \{V_i\})$ is a dynamical datum
over the abelian group $A=<g>$ and we compute the corresponding
dynamical twist. First let us prove the following technical
result.

\begin{lema}\label{tech-taft} Let $\eta \in \ku^{\times}$.
Define $\xi\in T(q)$ by
\begin{align}\label{xi1} \xi=1+ \sum_{i=1}^{n-1} \, a_i\, x^i g^{n-i},\end{align}
where $a_l=\eta^l \prod_{j=1}^l
\displaystyle{\frac{q^{n-j+1}}{q^j-1}} $ for $l=1\dots n-1$. Then
$\xi$ is invertible and  $\xi(g+\eta x)=g \xi$.
\end{lema}
\pf The proof that $\xi(g+\eta x)=g \xi$ is done by a
straightforward computation. Is easy to see that
$(\sum_{i=1}^{n-1} \, a_i\, x^i g^{n-i})^n=0$, this implies that
$\xi$ is invertible.

\epf

We shall denote $\xi_j=1+ \sum_{i=1}^{n-1} \, a_i\, x^i g^{n-i}$,
where $a_l=\displaystyle{\frac{1}{c^l q^{lj}}} \prod_{j=1}^l
\displaystyle{\frac{q^{n-j+1}}{q^j-1}} $.

\medbreak

Set $\chi_i$ the character of the group $A$ determined by
$\chi_i(g)=q^i$. So $$C(\chi_i)=\{\alpha\in T(q)^*: \langle\alpha,
gt\rangle= q^i \langle\alpha, t\rangle\;\, \text{ for all } t\in
T(q)\}.$$ Denote $V_{\chi_i}=V_i$ for all $i=0\dots n-1$.

\begin{prop} The collection $(\Ac(1,c), \{V_{\chi_i}\}_{i=0\dots n-1})$ is a
dynamical datum.
\end{prop}
\pf Since all representations $V_i$ are one-dimensional we can
identify the stabilizer $\St_{\Ac(1,c)}(V_i, V_j)$ with the set
$$D(i,j)=\{\alpha\in T(q)^*: q^{i-j}\,\langle\alpha, t\rangle=
\langle \alpha, (g+ \frac{1}{c q^j}x)t\rangle\;\, \text{ for all }
t\in T(q)\}.$$ Indeed by \cite[Lemma 2.8]{AM} $\alpha\ot
T\in\St_{\Ac(1,c)}(V_i, V_j)$ if and only if
$$\langle \alpha, k\_{-1}t\rangle \, T(k\_0\cdot v_i)=\langle\alpha,t \rangle
\, k\cdot T(v_i)$$ for all $k\in\Ac(1,c), t\in T(q)$. Since $V_i$
are one-dimensional we can assume that $ T(v_i)=v_j$, and taking
$k=y$ we get the result. Now we shall prove that there is a
$T(q)$-module isomorphism $D(i,j)\simeq C(\chi_{i-j}).$

Define $\omega_{i,j}: D(i,j)\to C(\chi_{i-j})$ by
$$\langle\omega_{i,j}(\alpha),t\rangle=\langle\alpha, \xi_j^{-1}
t\rangle.$$ The maps $\omega_{i,j}$ are well-defined, indeed if
$\alpha\in D(i,j)$ and $t\in T(q)$ then
\begin{align*} \langle\omega_{i,j}(\alpha), gt\rangle&= \langle\alpha
,\xi_j^{-1} gt\rangle=\langle\alpha , (g+ \frac{1}{c
q^j}x)\xi_j^{-1}t\rangle\\
&=q^{i-j}\,\langle\alpha, \xi_j^{-1}t\rangle=
q^{i-j}\,\langle\omega_{i,j}(\alpha), t\rangle.
\end{align*}
The second equality by Lemma \ref{tech-taft}. Thus
$\omega_{i,j}(\alpha)\in C(\chi_{i-j})$. Clearly $\omega_{i,j}$ is
a $T(q)$-module isomorphism. \epf

Now we compute the dynamical twist associated to $(\Ac(1,c),
\{V_{\chi_i}\})$. Let $X\in \Rep(T(q))$ and $x\in X[\chi_j]$,
$f\in X^*$ then
\begin{align*} \Psi(\chi_i, x)(f\ot
v_i)&=\omega^{-1}_{i,i+j}(\widetilde{f}^x)(1\ot\, v_i)=\langle
\widetilde{f}^x, \xi_{i+j}\rangle\; v_{i+j}\\
&=\langle f, \Ss(\xi_{i+j})\cdot x\rangle\; v_{i+j}.
\end{align*}
For any $i=0\dots n-1$ denote $P_i=P_{\chi_{i}}$. Let $ l,r=0\dots
n-1$ $\chi_l, \chi_r\in \widehat{A}$ and $f_1, f_2 \in T(q)^*$,
then $\Psi(\chi_i, I(\chi_i)(P_r\ot \, P_l))(f_1\ot f_2\ot\, v_i)$
is equal to
\begin{align*} \langle f_1, \Ss(\xi_{i+r+l})\_2I^2(\chi_i)P_l\rangle
\;\,\langle f_2, \Ss(\xi_{i+r+l})\_1I^1(\chi_i)P_r\rangle\;\,
v_{i+r+l} .
\end{align*}
On the other hand
$$\Psi(\chi_{i+r}, P_l)(\id\ot \Psi(\chi_{i}, P_r))(f_1\ot f_2\ot\, v_i)$$
is equal to
\begin{align*}\langle f_1, \Ss(\xi_{i+r+l})P_l\rangle
\;\,\langle f_2, \Ss(\xi_{i+r})P_r\rangle\;\, v_{i+r+l}.
\end{align*}
Hence
$$\Ss(\xi_{i+r+l})\_2 P_l\ot \, \Ss(\xi_{i+r+l})\_1P_r=\Ss(\xi_{i+r+l})
P_l I^{-2}(\chi_i)\ot \, \Ss(\xi_{i+r})P_r I^{-1}(\chi_i), $$ so
we deduce that
$$J(\chi_i)=\sum_{r,l} \, (\xi_{i-r-l})\_1 \xi^{-1}_{i-r-l} P_l \ot \,
 (\xi_{i-r-l})\_2  \xi^{-1}_{i-r} P_r. $$

\begin{rmk} It would be  interesting to prove that for each $d$ divisor of $n$
$(\Ac(d,c), \{V^d_i\})$, where $V^d_i$ are the irreducible
$\Ac(d,c)$-modules, is a dynamical datum. This result would
classify all dynamical twists since the categories
$_{\Ac(d,c)}\Mo$ are all exact indecomposable module categories
over $\Rep(T(q))$, \cite{eo}.
\end{rmk}

\end{document}